\documentclass[preprint,10pt]{elsarticle}
\usepackage{amsmath}
\usepackage{amsthm}
\usepackage{algorithm}
\usepackage{algpseudocode}

\newtheorem{theorem}{Theorem}[section]

\newtheorem{corollary}[theorem]{Corollary}

\newtheorem{example}[theorem]{Example}

\newproof{pf}{Proof}

\begin{document}
\begin{frontmatter}
\title{Infinity norm bounds for the inverse of
Nekrasov matrices using scaling matrices}
\author[rvt]{H. Orera\corref{cor1}}
\ead{hectororera@gmail.com}

\author[rvt]{J.~M. Pe\~{n}a}
\ead{jmpena@unizar.es}

\cortext[cor1]{Corresponding author}

\address[rvt]{Departamento de Matem\'{a}tica Aplicada/IUMA, Universidad de Zaragoza, Spain}

\begin{abstract}
For many applications, it is convenient to have good upper bounds for the norm of the inverse of a given matrix. In this paper, we obtain such bounds when $A$ is a Nekrasov matrix, by means of a scaling matrix transforming $A$ into a strictly diagonally dominant matrix. Numerical examples and comparisons with other bounds are included. The scaling matrices are also used to derive new error bounds for the linear complementarity problems when the involved matrix is  a Nekrasov matrix. These error bounds can improve considerably other previous bounds.

 
\end{abstract}
\end{frontmatter}

\section{Introduction}

Providing upper bounds for the infinity norm of the inverse of a  matrix has many potential applications in Computational Mathematics. For instance, for bounding the condition number of the matrix, for bounding errors in linear complementarity problems (cf. \cite{G-E}) or, in the class of $H$-matrices, for proving the convergence of matrix splitting and matrix multisplitting iteration methods for solving sparse linear systems of equations (cf. \cite{7}).

The class of Nekrasov matrices (see \cite{16} or Section 2) contains the class of strictly diagonally dominant matrices. Recent applications of Nekrasov matrices can be seen in  \cite{6}--\cite{SC}. Nekrasov matrices are $H$-matrices. Let us recall some related classes of matrices. 
A real matrix $A$ is a  nonsingular {\it $M$-matrix} if
its inverse is nonnegative and all its  off-diagonal entries are
nonpositive.  $M$-matrices form a very important class of matrices with applications to Numerical Analysis, Optimization, Economy, and Dynamic systems (cf. \cite{2}). Given a complex matrix $A=(a_{ij})_{1\leq i,j\leq n}$, its
{\it comparison matrix} $\mathcal{M}(A)=(\tilde{a}_{ij})_{1\leq i,j\leq
n}$ has entries $\tilde{a}_{ii}:=|a_{ii}|$ and
$\tilde{a}_{ij}:=-|a_{ij}|$ for all $j\ne i$ and $i,j=1,\ldots,n$.
We say that a complex matrix is an {\it $H$-matrix} if its comparison matrix
is a  nonsingular $M$-matrix. About a more general definition of $H$-matrix, see \cite{b}. A matrix
$A=(a_{ij})_{1\leq i,j\leq n}$ is {\it SDD} (strictly diagonally
dominant by rows) if $|a_{ii}|>\sum _{j\ne i}|a_{ij}|$ for all
$i=1,\ldots,n$. 

It is well--known that an SDD matrix is nonsingular and  that a square matrix $A$ is   an
 $H$-matrix if there exists a diagonal matrix $S$ with positive
diagonal entries such that $AS$ is SDD. The role of the scaling matrix is crucial, for instance, for the problem mentioned above of the convergence of iteration methods and also for the problem of eigenvalue localization (see \cite{SC}). This paper deals with the research of such scaling matrices $S$ for the particular case when $A$ is a Nekrasov matrix. The scaling matrix $S$ is applied to obtain upper bounds for the infinity norm of the inverse of a Nekrasov matrix. The paper is organized as follows.

 Section 2 constructs scaling matrices $S$ for Nekrasov matrices $A$, such that $AS$ is SDD. Section 3 applies the scaling matrices of Section 2 to derive upper bounds of $\left\Vert{A^{-1}}\right\Vert_\infty$, including an algorithm to obtain the corresponding bound.  Section 4 presents an improvement of the bound obtained in Section 3 and includes numerical examples, illustrating our bounds and comparing them with other previous bounds. We consider several test matrices previously considered in the literature, and we also consider some variants of these matrices. We also include a family of $3\times 3$ matrices showing that previous bounds can be arbitrarily large, in contrast to our bounds, which are always controlled.  Finally, we derive bounds for other norms. Section 5 illustrates the use of our  scaling matrices to derive new error bounds for the linear complementarity problems when the involved matrix is  a Nekrasov
 matrix. We avoid the restrictions of the bound in \cite{G-E} and we present a family of matrices for which our bound is a small
constant, in contrast to the bounds of \cite{c-x}, \cite{n1} and \cite{n2}, which can be arbitrarily
large.

We finish this introduction with some basic notations. Let $N:=\{1, \ldots,n\}$. Let $Q_{k,n}$ be the set of increasing sequence of $k$ positive integers in $N$. Given $\alpha,\beta\in Q_{k,n}$, we denote by $A[\alpha|\beta]$ the $k\times k$ submatrix of $A$ containing rows numbered by $\alpha$ and columns numbered by $\beta$. If $\alpha=\beta$, then we have the principal submatrix $A[\alpha]:=A[\alpha|\alpha]$. Finally, the diagonal matrix with diagonal entries $d_i$, $1\le i\le n$, will be denoted by ${\rm diag}(d_i)_{i=1}^n$.

\section{Scaling matrices}

Let us start by defining the concept of a Nekrasov matrix (see \cite{7,16,6}).  For this purpose, let us define recursively, for a complex matrix
$A=(a_{ij})_{1 \leq i,j \leq n}$ with $a_{ii}\neq 0$, for all $i=1,
\ldots,n$,
\begin{equation}\label{hs}
h_1(A):=\sum_{j\neq 1}|a_{1j}|,
\,\,\,h_i(A):=\sum_{j=1}^{i-1}|a_{ij}|{h_j(A) \over
|a_{jj}|}+\sum_{j=i+1}^n|a_{ij}|, \quad i=2, \ldots,n.
\end{equation}
 We say that $A$ is a {\it Nekrasov matrix} if
$|a_{ii}|>h_i(A)$ for all $i\in N$. 
It is well--known that a Nekrasov
matrix is a nonsingular $H$-matrix \cite{16}.  So, there exists a positive
diagonal matrix $S$ such that $AS$ is {\it SDD}. In particular, Nekrasov matrices can be characterized in terms of these  scaling matrices (see Theorem 2.2 of \cite{K2}). Once we have found a scaling matrix, we can use it to derive infinity norm bounds for the inverse of Nekrasov matrices, which may be useful for many problems, as recalled in the Introduction. In fact, the problem of bounding the infinity norm of the inverse of a  Nekrasov matrix  has attracted great attention recently (see \cite{7,G,K,K3,L}).\\


 In this section we are introducing two methods that allow us to build a scaling matrix for any given Nekrasov matrix.

\begin{theorem}\label{escaladoSDD}
Let $A=(a_{ij})_{1\leq i,j \leq n}$ be a Nekrasov matrix. Then the matrix
$$S=\begin{pmatrix}
\frac{h_1(A) + \epsilon_1}{|a_{11}|} &  & \\ 
 &\ddots  & \\ 
 &  & \frac{h_n(A) + \epsilon_n}{|a_{nn}|}
\end{pmatrix}, $$
\begin{flalign*}
 \text{with } 
\left\{\begin{array}{ll}
\epsilon_1>0,\\
0<\epsilon_i\leq|a_{ii}| - h_i(A), & \epsilon_i>\sum_{j=1}^{i-1}\frac{|a_{ij}|\epsilon_j}{|a_{jj}|} \ \ \text{  for } i=2,\ldots,n,\
\end{array}\right.  &&  
\end{flalign*}
is a positive diagonal matrix  such that $AS$ is SDD.
\end{theorem}

\begin{proof}
Let us start by proving that there exist $\epsilon_1,\ldots,\epsilon_n$ satisfying the conditions stated above. We consider, as a first choice, $\epsilon_i:=|a_{ii}|-h_i(A)$ for $i=1,\ldots,n$. 
If  $\epsilon_i>\sum_{j=1}^{i-1}\frac{|a_{ij}|\epsilon_j}{|a_{jj}|}$ for all $i=2,\ldots,n$ we have finished. Otherwhise, let $i>1$ be the first index such that the inequality does not hold. Then we substitute $\epsilon_j$ by $\frac{\epsilon_j}{\hat{M}_i}$, with $j=1,\ldots,i-1$, where $\hat{M}_i$ is a positive number such that the inequality is satisfied. The inequalities checked at earlier steps  remain true. We continue this process  until the inequality  holds for all $i=2,\ldots,n$.

The diagonal matrix $S$ is positive because $h_i(A)\geq 0$ and $\epsilon_i>0$. The entry $(i,j)$ of $AS$ is $a_{ij}\frac{h_j(A) + \epsilon_j}{|a_{jj}|} $. In order to prove that  $AS$ is SDD we start by checking that the condition is true for the  $n$th row:
$$\sum_{j=1}^{n-1}|a_{nj}|\frac{h_j(A)+\epsilon_j}{|a_{jj}|}=\underbrace{\sum_{j=1}^{n-1}|a_{nj}|\frac{h_j(A)}{|a_{jj}|}}_{h_n(A)} +\sum_{j=1}^{n-1}|a_{nj}|\frac{\epsilon_j}{|a_{jj}|}<h_n(A)+\epsilon_n=|(AS)[n]| $$
The condition holds for the  row $n-1$:
\begin{equation}
\begin{split}
\nonumber
\sum_{j=1}^{n-2}|a_{n-1,j}|\frac{h_j(A)+\epsilon_j}{|a_{jj}|} + |a_{n-1, n}|\underbrace{\frac{h_n(A)+\epsilon_n}{|a_{nn}|}}_{\leq 1} &\leq h_{n-1}(A) + \sum_{j=1}^{n-2}|a_{n-1 ,j}|\frac{\epsilon_j}{a_{jj}}\\
<h_{n-1}(A) + \epsilon_{n-1}=|(AS)[n-1]| &
\end{split}
\end{equation}
The first inequality is due to the hypothesis $\epsilon_n\leq |a_{nn}| - h_n(A) $, which implies $\frac{h_n(A) + \epsilon_n}{|a_{nn}|}\leq 1$.
In general, considering the $i$th row for $2\leq i <n-1$:
\begin{equation}
\begin{split}
\nonumber
\sum_{j=1}^{i-1}|a_{ij}|\frac{h_j(A)+\epsilon_j}{|a_{jj}|} + \sum_{j=i+1}^{n}|a_{ij}|\frac{h_j(A)+\epsilon_j}{|a_{jj}|} &\leq h_i(A) + \sum_{j=1}^{i-1}|a_{ij}|\frac{\epsilon_j}{|a_{jj}|}\\
<h_i(A)+\epsilon_i = |(AS)[i]| &
\end{split}
\end{equation}
and, when $i=1$:
$$\sum_{j=2}^{n}|a_{1j}|\frac{h_j(A)+\epsilon_j}{|a_{jj}|} \leq h_1(A) <h_1(A) + \epsilon_1=|(AS)[1]|. $$
The inequality for the $i$th row is proven using that $\epsilon_j\leq |a_{jj}|-h_j(A)$ for $j=i+1,\ldots,n$ and $\epsilon_i>\sum_{j=1}^{i-1}\frac{|a_{ij}|\epsilon_j}{|a_{jj}|}$. If $i=1$, the last inequality is reduced to $\epsilon_1>0$ .
\end{proof}

In  
Theorem \ref{escaladoSDD} we introduced a diagonal matrix $S$ that transforms any Nekrasov matrix into an SDD matrix. Its construction  implied the search of the parameters $\epsilon_i$ for $i\in N$. Taking into account the existence of nonzero entries in the upper triangular part of  a Nekrasov matrix, we can build a new scaling matrix $S$, simpler in many cases, whose product with a Nekrasov matrix is also SDD. 

 \begin{theorem}\label{corolarioSDD} Let $A=(a_{ij})_{1\leq i,j \leq n}$ be a Nekrasov matrix, and let $k\in N$ be the first index such that there does not exist $j>k$  with $ a_{kj}\neq 0$. Then the matrix $S={\rm diag}(s_i)_{i=1}^n$ with $s_i:=\frac{h_i(A) + \epsilon_i}{|a_{ii}|}$ and with
$$\left\{\begin{array}{ll}
\epsilon_i=0, & i=1,\ldots,k-1,\\
0<\epsilon_i<|a_{ii}| - h_i(A), & \epsilon_i>\sum_{j=k}^{i-1}\frac{|a_{ij}|\epsilon_j}{|a_{jj}|} \ \ \text{  for } i=k,\ldots,n,\
\end{array}\right. $$
is a positive diagonal matrix such that $AS$ is SDD. 
\end{theorem}
\begin{proof}
Let us start by showing that there exist $\epsilon_1,\ldots,\epsilon_n$ satisfying the conditions stated above. Since $A$ is a Nekrasov matrix, we have that $|a_{ii}|>h_i(A)$ for $i=1,\ldots,n$. The existence of $\epsilon_1=\ldots=\epsilon_{k-1}=0$ is trivial and, following the constructive proof of the existence of these parameters given in Theorem \ref{escaladoSDD}, we can deduce the existence of $\epsilon_k,\ldots,\epsilon_n$. It remains to prove that $AS$ is an SDD matrix, which can be done analogously to the proof of Theorem \ref{escaladoSDD}.

Let us first consider  the $i$th row, when $i<k$. Since $i<k$, there exists an entry $a_{ij}\neq 0$ with $i<j$. Taking also into account that $\epsilon_j=0$ for all $j<i$ and that $h_j(A)+\epsilon_j < |a_{jj}|$ for all $j=i+1,\ldots,n,$ we deduce that:

\begin{equation}
\nonumber
\begin{split}
\sum_{j=1}^{i-1}|a_{ij}|\frac{h_j(A)+\epsilon_j}{|a_{jj}|} + \sum_{j=i+1}^{n}|a_{ij}|\frac{h_j(A)+\epsilon_j}{|a_{jj}|}\\
= \sum_{j=1}^{i-1}|a_{ij}|\frac{h_j(A)}{|a_{jj}|} + \sum_{j=i+1}^{n}|a_{ij}|\frac{h_j(A)+\epsilon_j}{|a_{jj}|}\\
 < \sum_{j=1}^{i-1}|a_{ij}|\frac{h_j(A)}{|a_{jj}|} + \sum_{j=i+1}^{n}|a_{ij}|=h_i(A)=|(AS)[i]|.
\end{split}
\end{equation} 
For the $k$th row we have that $a_{kj}=0$ for every $j>k$ and so:
\begin{equation}
\nonumber
\begin{split}
\sum_{j=1}^{k-1}|a_{kj}|\frac{h_j(A)+\epsilon_j}{|a_{jj}|} &+ \sum_{j=k+1}^{n}|a_{kj}|\frac{h_j(A)+\epsilon_j}{|a_{jj}|}=\sum_{j=1}^{k-1}|a_{kj}|\frac{h_j(A)}{|a_{jj}|} \\
  = h_k(A)<h_k(A)&+\epsilon_k=|(AS)[k]|.
\end{split}
\end{equation}
It just remains to check the $i$th rows, when $i>k$. Since $h_j(A)+\epsilon_j<|a_{jj}|$ for all $j>i(>k)$, we have, by the choice of $\epsilon_i$:

\begin{equation}
\nonumber
\begin{split}
\sum_{j=1}^{i-1}|a_{ij}|\frac{h_j(A)+\epsilon_j}{|a_{jj}|} + \sum_{j=i+1}^{n}|a_{ij}|\frac{h_j(A)+\epsilon_j}{|a_{jj}|} &\leq h_i(A) + \sum_{j=k}^{i-1}|a_{ij}|\frac{\epsilon_j}{|a_{jj}|}\\
<h_i(A) + \epsilon_i = |(AS)[i]|.  
\end{split}
\end{equation}

\end{proof}

The particular case $k=n$ corresponds to  a diagonal matrix with $\epsilon_1=\ldots=\epsilon_{n-1}=0$ and $\epsilon_{n}\in \left(0, |a_{nn}|-h_n(A)\right)$. This scaling matrix was already introduced in \cite{G-E} and it was used to derive an error bound for linear complementarity problems of Nekrasov matrices. In the following section, we shall apply the scaling matrices derived in this section to the problem of bounding the norm of the inverse of a Nekrasov matrix.

\section{Bounding $\left\Vert{A^{-1}}\right\Vert_\infty$}

With an adequate scaling matrix $S$ (given by Theorem \ref{escaladoSDD} or Theorem \ref{corolarioSDD}) we can obtain the desired bound for the inverse of a Nekrasov matrix $A$ considering the product $AS$. For this purpose, we are going to use the following result introduced by J. M. Varah in \cite{Varah}:

\begin{theorem}\label{Varah}
If $A$ is SDD and  $\alpha:=\min_k (|a_{kk}|-\sum_{j\neq k}|a_{kj}|), \ then\linebreak \left\Vert{A^{-1}}\right\Vert_\infty<1/\alpha.$
\end{theorem}

Theorem \ref{Varah} gives a bound for the infinity norm of the inverse of an SDD matrix. This theorem, jointly with  the scaling matrices introduced in Section 2, allows us to deduce Theorem \ref{lacota}.

\begin{theorem}\label{lacota} Let $A=(a_{ij})_{1\leq i,j \leq n}$ be a Nekrasov matrix. Then
\begin{equation}
\label{cotanek}
 \left\Vert{A^{-1}}\right\Vert_\infty \leq  \frac{\max_{i\in N}\left(\frac{h_i(A)+\epsilon_i}{|a_{ii}|}\right)}{\min_{i\in N}(\epsilon_i-w_i+p_i)}, 
\end{equation}
where $(\epsilon_1,\ldots,\epsilon_n)$ are given by Theorem \ref{escaladoSDD} or Theorem \ref{corolarioSDD}, $w_i:=\sum_{j=1}^{i-1}|a_{ij}|\frac{\epsilon_j}{|a_{jj}|}$,  and $p_i:=\sum_{j=i+1}^{n}|a_{ij}|\frac{|a_{jj}|-h_j(A)-\epsilon_j}{|a_{jj}|}$ for all $i\in N$.
\end{theorem}

\begin{proof}
We choose a diagonal matrix $S$ following either Theorem \ref{escaladoSDD} or Theorem \ref{corolarioSDD} and we deduce the following inequality:
\begin{equation}\label{dem.cota} 
\left\Vert A^{-1}\right\Vert_\infty = \left\Vert{S(S^{-1}A^{-1})}\right\Vert_\infty = \left\Vert{S(AS)^{-1}}\right\Vert_\infty \leq \left\Vert{S}\right\Vert_\infty \left\Vert{(AS)^{-1}}\right\Vert_\infty .
\end{equation}
The matrix $S$ is diagonal, so its infinity norm is given by $\max_{i\in N}\left(\frac{h_i(A)+\epsilon_i}{|a_{ii}|}\right)$. Since $AS$ is SDD  we can apply Theorem \ref{Varah} to $\left\Vert{(AS)^{-1}}\right\Vert_\infty$. For this purpose, we need to compute for each $i=1,\ldots,n$:
 \begin{equation}
 \nonumber
 \begin{split}
 h_i(A)+ \epsilon_i -\sum_{j\neq i}|a_{ij}|\frac{h_j(A)+\epsilon_j}{|a_{jj}|} & = \epsilon_i - \sum_{j=1}^{i-1}|a_{ij}|\frac{\epsilon_j}{|a_{jj}|} + \sum_{j=i+1}^{n}|a_{ij}|\frac{|a_{jj}|-h_j(A)-\epsilon_j}{|a_{jj}|}\\
 &=\epsilon_i-w_i+p_i,\\
 \end{split}
  \end{equation}
 where we have substituted  $h_i(A)$ by the expression given by \eqref{hs}. 
\end{proof}

Since the diagonal matrix $S$ satisfies $\left\Vert{ S }\right\Vert_\infty\leq 1$, we can substitute the numerator of the bound (\ref{cotanek}) by one and obtain the following result:

\begin{corollary} \label{cf}
 Let $A=(a_{ij})_{1\leq i,j \leq n}$ be a Nekrasov matrix. Then
\begin{equation}
\nonumber
 \left\Vert{A^{-1}}\right\Vert_\infty \leq  \frac{1}{\min_{i\in N}(\epsilon_i-w_i+p_i)}, 
\end{equation}
where $(\epsilon_1,\ldots,\epsilon_n)$ are given by Theorem \ref{escaladoSDD} or Theorem \ref{corolarioSDD}, $w_i:=\sum_{j=1}^{i-1}|a_{ij}|\frac{\epsilon_j}{|a_{jj}|}$,  and $p_i:=\sum_{j=i+1}^{n}|a_{ij}|\frac{|a_{jj}|-h_j(A)-\epsilon_j}{|a_{jj}|}$ for $i\in N.$
\end{corollary}

In Table \ref{CompCost} we present the computational cost of the bound \eqref{cotanek} using the matrix $S$ given by Theorem \ref{corolarioSDD}. The cost depends on the index $k$. Two extreme cases are studied separately. The first one, $k=n$, corresponds to the simplest case, where $\epsilon_i=0$ for $i=1,\ldots,n-1$. The second one corresponds to $k=1$ and it uses a diagonal matrix $S$ with $\epsilon_i\neq 0$ for all $i\in N$. In fact, in this case the  diagonal matrix $S$  also satisfies the definition given by Theorem \ref{escaladoSDD}.

\begin{table}[h]
\centering
\begin{tabular}{|c|c|c|c|}
 \hline
Operations &general &$k=n$ &$k=1$  \\ 
 \hline
additions/subtractions & $ \frac{3n^2+n+2}{2}+\frac{(n-k-1)(n-k)}{2}$  &  $ \frac{3n^2+n+2}{2}$ &$ 2n^2-n + 2$\\
multiplications & $\frac{7n^2 + 9n + 4}{2} + \frac{5k^2 - 10kn - 11k}{2}$  & $n(n-1) $  & $\frac{7n^2 -n -2}{2}$\\
quotients & $2n-1 + 2(n-k)$  & $2n-1$& $4n-3$\\
 \hline
\end{tabular}
 \caption{Computational cost of \eqref{cotanek}}
 \label{CompCost}
\end{table}

The particular cases $k=n$ and $k=1$ have the lowest and biggest computational cost, respectively. Table \ref{leadingTerm} shows
the leading term $T$ of the computational cost in these cases.\\

\begin{table}[h]
\centering
\begin{tabular}{|c|c|c|c|}
 \hline
  &$k=n$ &$k=1$  \\ 
 \hline
$T$ &    $ \frac{5}{2}n^2 $ &$ \frac{11}{2}n^2 $\\
 \hline
\end{tabular}
\caption{Leading term of the computational cost of \eqref{cotanek}}
\label{leadingTerm}
\end{table}

\begin{algorithm}
\caption{nektoSDD - Computing bound (\ref{cotanek})}\label{nektoSDD}
\begin{algorithmic}
\State \bf{Input:}  $A=(a_{ij})_{1\leq i,j \leq n}$,  $t$ 
\For{$i=1:n$}
	\State $ h_{i}=\sum _{j=1}^{i-1}|a_{ij}|k_j$ 
	\State $r= \sum _{j=i+1}^{n}|a_{ij}|$
	\If{$r==0 , J==0$} \Comment{\normalfont Find the first row such that $\epsilon_i>0$}
	    \State J=i; 
	\EndIf
	\State $ h_i=h_i+r$
	\State $\Delta_i =|a_{ii}| - h_i$
	\State $ k_{i}=h_i/|a_{ii}|$ 
\EndFor
\State $\epsilon_K=t\Delta_K$
\State $w_1=\ldots=w_K=0$ \Comment{\normalfont If $i\leq K$, we have that $w_i=0$}
\For{$i=K+1:n$}
    \State $\epsilon_i=t\Delta_i$
    \State $p_j=\epsilon_j /|a_{jj}|$
    \State $w_i=\sum_{j=K}^{i-1}|a_{ij}|p_j$
    \If{$w_i-\epsilon_i>0$}
        \State $M=1/2w_i$
        \For{$j=K:i-1$}
            \State $\epsilon_j=\epsilon_j\epsilon_i M$
            \State $w_j=w_j\epsilon_i M$
        \EndFor
        \State $w_i=\epsilon_i/2$
    \EndIf
\EndFor 
\For{i=n:-1:2}
    \State $S_i=\epsilon_i-w_i + \sum_{j=i+1}^n |a_{ij}|f_j $ 
    \State $f_i=(\Delta_i - \epsilon_i)/|a_{ii}|$
\EndFor
\State $S_1=\epsilon_1-w_1 + \sum_{j=2}^n |a_{1j}|f_j $ 
\State $\displaystyle Bound=\frac{\max_{i\in N}{\{(k_i+\epsilon_i/|a_{ii}|\}}}{\min_{i\in N}{\{S_i\}}}$
\end{algorithmic}
\end{algorithm}
Now we are going to introduce Algorithm \ref{nektoSDD}, which allows us to compute the bound (\ref{cotanek}) choosing $\epsilon_i$ with $i\in N$ following Theorem \ref{corolarioSDD}. It corresponds to  Theorem \ref{escaladoSDD} when $k=1$. In order to compute this bound, the algorithm needs to give some initial values to $\epsilon_1,\ldots,\epsilon_n$.  These parameters are initialized
 with either $0$ or $t(|a_{ii}|-h_i(A))$, where $t \in (0,1)$. It could be useful to choose a different  scalar $t$ for each $\epsilon_i$. However, it is not clear how to choose  their values, and for many matrices, such as those included in Section \ref{sec_tests}, we have that   $\epsilon_i=0$ for $i=1,\ldots,n-1$. In this case, we also consider in Section \ref{sec_tests} the possibility of choosing $\epsilon_n$ as the middle point of its interval, that is, $\epsilon_n =\frac{\Delta_n}{2}$, where $\Delta_n = |a_{nn}| - h_n(A)$.

\section{Improvements, numerical tests and bounds for other norms}
\label{sec_tests}

In the previous section, we derived the bound \eqref{cotanek} for the infinity norm of the inverse of a Nekrasov matrix $A$. For this purpose, we first obtained an adequate scaling matrix $S$ and then we applied the well--known Varah's bound of Theorem \ref{Varah} to the matrix $AS$. Nevertheless, any bound applicable to SDD matrices could be applied to $AS$, and a different choice would lead us to a different bound. In order to illustrate this fact, we are also going to use the bound introduced in \cite{K} for Nekrasov matrices, which in particular improves Varah's bound for SDD matrices (as proven in Theorem 2.4 of \cite{K}):

\begin{equation}
\label{cotaK}
 \left\Vert{A^{-1}}\right\Vert_\infty \leq  \max_{i\in N}\frac{z_i(A)}{|a_{ii}|-h_i(A)}, 
\end{equation}
\begin{equation*}
z_1(A):=1,
\,\,\,z_i(A):=\sum_{j=1}^{i-1}|a_{ij}|{z_j(A) \over
|a_{jj}|} + 1, \quad i=2, \ldots,n.
\end{equation*}
As in \eqref{dem.cota}, the new bound for $\left\Vert{A^{-1}}\right\Vert_\infty$ reduces to the product of $\left\Vert{S}\right\Vert_\infty$ and the bound to $ \left\Vert{(AS)^{-1}}\right\Vert_\infty$ obtained by \eqref{cotaK}. In fact, taking into account that $z_i(AS)=z_i(A)$ for all $i\in N$, the explicit form of this new bound is:

\begin{equation}\label{cotarev}
 \left\Vert{A^{-1}}\right\Vert_\infty\leq \max_{i\in N }\left(\frac{h_i(A)+\epsilon_i}{|a_{ii}|}\right)\max_{i\in N}\frac{ z_i(A)}{(h_i(A)+\epsilon_i - h_i(AS)) }.
\end{equation}
 As  shown by the following numerical experiments, this change gives a better  bound  whenever $S$ follows Theorem \ref{corolarioSDD}. However, in general the substitution of Varah's bound  is going to  increase the computational cost of the bound, while the bound \eqref{cotanek} using Theorem \ref{escaladoSDD} is  not significantly improved. Analogously to \eqref{cotarev}, if better bounds than \eqref{cotaK} for SDD matrices are obtained, then they can be also combined with our bound of Theorem \ref{corolarioSDD} to derive sharper bounds than \eqref{cotarev}, although the computational cost can increase again. 

Recent articles have studied the problem of finding bounds for the infinity norm of the inverse of a Nekrasov matrix. In \cite{7}, two bounds  are introduced and tested with the following six matrices:

\begin{center}
$\begin{array}{ccc}
 A_1=\begin{pmatrix}
-7 &1 &-0.2& 2\\
7 & 88 & 2 &-3\\
2 & 0.5 & 13 & -2\\
0.5 & 3 & 1 & 6
\end{pmatrix},
&
A_2=\begin{pmatrix}
8 &1 &-0.2 &3.3\\
7 &13 &2 &-3\\
-1.3& 6.7& 13 &-2\\
0.5 &3& 1 &6
\end{pmatrix},\\
\\
A_3=\begin{pmatrix}
21 & -9.1 &-4.2 &-2.1\\
-0.7& 9.1 &-4.2 &-2.1\\
-0.7 &-0.7 &4.9 &-2.1\\
-0.7 &-0.7& -0.7& 2.8
\end{pmatrix}, & A_4=\begin{pmatrix}
5 & 1& 0.2 &2\\
1 &21 &1& -3\\
2 &0.5 &6.4 &-2\\
0.5& -1& 1 &9
\end{pmatrix}, \\
\\
A_5=\begin{pmatrix}
6 &-3 &-2\\
-1& 11& -8\\
-7& -3 &10
\end{pmatrix},&A_6=\begin{pmatrix}
8 &-0.5& -0.5 &-0.5\\
-9& 16 &-5 &-5\\
-6 &-4& 15 &-3\\
-4.9 &-0.9 &-0.9& 6
\end{pmatrix}.
\end{array}$\
\end{center}

In more recent works, such as \cite{G,K,L}, improvements of these bounds are developed and  tested using also these matrices. Since the scaling matrices introduced in Section 2 allowed us to derive different bounds, we are going to compare them with the  results obtained in some of the mentioned  papers.

 We have included results from \cite{G,K}. The bound \eqref{cotaK} (which corresponds to the bound 2.4 of \cite{K})   improves those obtained in \cite{7} for Nekrasov matrices (as proven in Theorem 2.3 of \cite{K}). Theorem 9 of \cite{G} gives a sharper bound in some cases, and so we also include it in our comparison.

 Table \ref{comp} gathers the different bounds. The first row shows the exact infinity norm of the matrices. The data included in the second (corresponding to bound \eqref{cotaK}) and third rows are borrowed from the  articles that achieved the sharpest bounds. The other rows contain our results, obtained with bounds \eqref{cotanek} and \eqref{cotarev}.  In the last case $S$ was given by Theorem \ref{escaladoSDD} while in the other cases  the diagonal matrix $S$  followed Theorem \ref{corolarioSDD}. Excluding the case where $\epsilon_n=\Delta_n/2$, our bounds  used an appropriate choice of parameters.

\begin{table}[h]
\centering
\begin{tabular}{|c|c|c|c|c|c|c|}
 \hline
Matrix & $A_1$  & $A_2$& $A_3$& $A_4$ & $ A_5$ & $ A_6$ \\ 
 \hline
Exact norm &0.1921 &0.2390  & 0.8759 & 0.2707 & 1.1519  &0.4474 \\
\eqref{cotaK} &0.2632 &0.5365  & 0.9676 & 0.5556 & 1.4138  &0.4928 \\
Theorem 9 of \cite{G} & 0.2505  &  0.5365 & 0.9676 & 0.5038 & 1.4138 & 0.4928 \\
\eqref{cotanek}, $\epsilon_n=\Delta_n/2$ & 0.6398 &1.4406  & 1.5527 & 0.7264 & 1.2974  & 1.2893  \\
\eqref{cotarev}, $\epsilon_n=\Delta_n/2$ & 0.4992 &0.7422  & 1.0632 & 0.5596& 1.2809  & 1.2893  \\
\eqref{cotanek}, Theorem \ref{corolarioSDD} &0.3474 &0.8894  & 1.3325 & 0.4484 & 1.1658  &1.0796 \\
\eqref{cotarev}, Theorem 2.2 &0.3074 &0.5684 & 0.9735 & 0.3817 & 1.1658 & 1.0436 \\
\eqref{cotanek}, Theorem \ref{escaladoSDD} &0.2354 &0.5260  & 0.9273 & 0.3168& 1.1588  &0.4527 \\
 \hline
\end{tabular}
\caption{Upper bounds of $\vert\vert{A^{-1}}\vert\vert_\infty$}
\label{comp}
\end{table}

Looking at the rows corresponding to Theorem \ref{corolarioSDD} we can observe that the obtained bounds are better for $A_4$ and $A_5$, but they are worse in the other cases. With the choice of a diagonal matrix $S$ following Theorem \ref{escaladoSDD} and bound \eqref{cotanek} we  obtained a better bound for every matrix. This option seems superior to the other possibilities. However, it has an intrinsic problem: the choice of the parameters $\epsilon_i$ for $i=1,\ldots,n$. Given the right parameters, the obtained bound  is excellent. But a bad choice of these values may give a useless bound. In general, it is not clear how to find the optimal values using Theorem \ref{escaladoSDD}.\\

Performing more numerical tests, we have seen that the bounds introduced in this paper may be particularly useful when the considered Nekrasov matrix is far from satisfying $|a_{ii}|>h_i(A)$ for some $i\in N\backslash \{n\}$. Looking at the bound for SDD matrices introduced by Varah (Theorem \ref{Varah}), we can observe that it depends on all the row sums of the comparison matrix. In particular, bounds for Nekrasov matrices based on Varah's bound seem to be inversely proportional to $|a_{ii}|-h_i(A)$ for some indices $i\in N$. In order to illustrate this fact, we have modified one entry of all previous examples and we present  the bounds obtained for the inverses of these new matrices in Table \ref{alt}.

\begin{center}
$\begin{array}{ccc}
 \hat{A}_1=\begin{pmatrix}
-7 &1 &\bf{-3.9}& 2\\
7 & 88 & 2 &-3\\
2 & 0.5 & 13 & -2\\
0.5 & 3 & 1 & 6
\end{pmatrix},
&
\hat{A}_2=\begin{pmatrix}
8 &1 &-0.2 &3.3\\
7 &13 &2 &-3\\
\bf{-11}& 6.7& 13 &-2\\
0.5 &3& 1 &6
\end{pmatrix},\\
\\
\hat{A}_3=\begin{pmatrix}
21 & -9.1 &-4.2 &-2.1\\
-0.7& 9.1 &-4.2 &\bf{-4.2}\\
-0.7 &-0.7 &4.9 &-2.1\\
-0.7 &-0.7& -0.7& 2.8
\end{pmatrix}, & \hat{A}_4=\begin{pmatrix}
5 & 1& 0.2 &2\\
1 &21 &1& -3\\
2 &0.5 &6.4 &-2\\
0.5& -1& \bf{15} &9
\end{pmatrix}, \\
\\
\hat{A}_5=\begin{pmatrix}
6 &-3 &-2\\
-1& \bf{9}& -8\\
-7& -3 &10
\end{pmatrix},&\hat{A}_6=\begin{pmatrix}
8 &-0.5& -0.5 &-0.5\\
\bf{-31.9}& 16 &-5 &-5\\
-6 &-4& 15 &-3\\
-4.9 &-0.9 &-0.9& 6
\end{pmatrix}.
\end{array}$\
\end{center}

\begin{table}[h]
\centering
\begin{tabular}{|c|c|c|c|c|c|c|}
 \hline
Matrix & $\hat{A}_1$  & $\hat{A}_2$& $\hat{A}_3$& $\hat{A}_4$ & $ \hat{A}_5$ & $ \hat{A}_6$ \\ 
 \hline
Exact norm &0.2385 &0.9827   & 1.0997   & 0.2848  & 2.4545  &  0.9144 \\
\eqref{cotaK} &10.0000  &16.2005   &   5.5357 &  8.7889  & 7.0000  & 266.0000  \\
Theorem 9 of \cite{G} & 0.3979   &  16.2005  &   5.5357 &  8.7889  & 7.0000 & 266.0000   \\
\eqref{cotanek}, $\epsilon_n=\Delta_n/2$ & 1.2345 & 2.2098 &  2.3120  & 17.0569   & 5.5208 & 2.6020 \\
\eqref{cotarev}, $\epsilon_n=\Delta_n/2$ & 0.6144 & 1.2071 &  1.6377  &  3.1074  & 5.5208 & 2.6020 \\
\eqref{cotanek}, Theorem \ref{corolarioSDD} &0.8230 & 1.4732 & 2.1018 & 10.2316 &4.1085   & 2.0316\\
\eqref{cotarev}, Theorem 2.2 &0.5344 & 0.9923& 1.5203 & 3.0603 &3.4717  & 1.9119\\
\eqref{cotanek}, Theorem \ref{escaladoSDD} &0.3262 &  1.2642 & 1.1479 & 6.6456 & 2.6180  &2.0316 \\

 \hline
\end{tabular}
\caption{Upper bounds of $\vert\vert{A^{-1}}\vert\vert_\infty$}
\label{alt}
\end{table}    
    
 In Table \ref{alt} we observe that bound \eqref{cotanek} is lower than  \eqref{cotaK} and the bound of \cite{G} even   with the choice of $\epsilon_n$ as the middle point  for matrices $\hat{A}_2$,$\hat{A}_3$,$\hat{A}_5$ and $\hat{A}_6$.  We also obtained tight bounds for the norm of the inverse of $\hat{A}_1$. The remaining case, $\hat{A}_4$, was built increasing significantly an entry of the last row. As a consequence,   all  bounds compared in Table \ref{alt} obtained weaker results than in Table \ref{comp}. For $\hat{A}_6$, bounds of \cite{G}  and \cite{K} (corresponding to bound \eqref{cotaK}) are very high while our bounds (using \eqref{cotanek}) are all controlled. This phenomenon will be also illustrated with the following family of $3\times 3$ matrices.\\
 
 \begin{example} Let us consider the family of matrices
 \begin{equation}
 \label{ex}
 A=\begin{pmatrix}
 4 & 2 & 1\\ 
 \frac{4}{3}-\varepsilon & 2 & 1\\ 
 1 & 1 & 2
 \end{pmatrix},
 \end{equation}
 where $0<\varepsilon <\frac{1}{10}$. In this case, $\vert\vert{A^{-1}}\vert\vert_\infty< 1.4167$, $h_1(A)=3$, $h_2(A)=2-\frac{3}{4}\varepsilon$ and $h_3(A)=\frac{7}{4}-\frac{3}{8}\varepsilon$.
 Then the bounds (2.4) of \cite{K} and Theorem 9 of \cite{G} coincide and are equal to $\frac{16}{9 \varepsilon} - \frac{1}{3}$. We can observe that this  bound is arbitrarily large when $\varepsilon \to 0$. However, our bounds remain controlled. In fact, \eqref{cotanek} with $\epsilon_3=\Delta_3/2$ (and $\epsilon _1=0$, $\epsilon _2=0$) gives the bound $16\left(\frac{1-(3\varepsilon/8)}{1+(3\varepsilon/2)}\right)$, \eqref{cotanek} in Theorem \ref{corolarioSDD} is equal to  $12\left(\frac{1-(3\varepsilon/8)}{1+(3\varepsilon/2)}\right)$ and \eqref{cotanek} in Theorem \ref{escaladoSDD} can become equal to 12.
  \end{example}
  
  We finish this section by applying Theorem \ref{lacota} to derive bounds for other norms. The first result is obtained from applying Theorem \ref{lacota}  to the transpose matrix.
  
  \begin{corollary} \label{cf1} 
  Let $A=(a_{ij})_{1\leq i,j \leq n}$ be a  matrix with $A^T$ Nekrasov. Then
\begin{equation}
\nonumber
 \left\Vert{A^{-1}}\right\Vert_1\leq \frac{\max_{i\in N}\left(\frac{h_i(A^T)+\bar\epsilon_i}{|a_{ii}|}\right)}{\min_{i\in N}(\bar \epsilon_i-\bar w_i+\bar p_i)}, 
\end{equation}
where $\bar \epsilon_i,\bar w_i, \bar p_i$ are  the parameters   $\epsilon_i, w_i,  p_i$    of   theorems \ref{corolarioSDD} and \ref{lacota} corresponding to $A^T$.
\end{corollary}

\begin{corollary} \label{cf2} 
  Let $A=(a_{ij})_{1\leq i,j \leq n}$ be a  matrix with $A$ and $A^T$ Nekrasov and let $\sigma _n(A)$ be its minimal singular value. Then
  \begin{equation}
\label{sigma}
 \sigma_n(A)=\left\Vert{A^{-1}}\right\Vert_2^{-1} \ge    \sqrt{\frac{\min_{i\in N}( \epsilon_i- w_i+ p_i)\min_{i\in N}(\bar \epsilon_i-\bar w_i+\bar p_i)}{\max_{i\in N}\left(\frac{h_i(A)+\epsilon_i}{|a_{ii}|}\right)\max_{i\in N}\left(\frac{h_i(A^T)+\bar\epsilon_i}{|a_{ii}|}\right)}}, 
\end{equation}
where $ \epsilon_i, w_i,  p_i,\bar \epsilon_i,\bar w_i, \bar p_i$ are given  in theorems \ref{corolarioSDD} and \ref{lacota} and Corollary \ref{cf1} .
\end{corollary}

\begin{proof}It is a consequence of the well--known facts that, for a nonsingular matrix $M$, its  minimal singular value coincides with $\left\Vert{M^{-1}}\right\Vert_2^{-1}$ and that $\Vert M\Vert_2^{2}\le \Vert M\Vert_1\Vert M\Vert_\infty $.
\end{proof}
  In the following example we apply  Corollary \ref{cf2} to the suitable matrices from the previous experiments, $A_3$ and $A_4$. The matrix $A_3$ is an SDD matrix whose transpose is Nekrasov, while $A_4$ is an SDD matrix whose transpose is also SDD. 
  
  \begin{example}  Corollary \ref{cf2} gives a lower bound for the minimal singular value of a Nekrasov matrix whose transpose is also a Nekrasov matrix. We can apply this result to $A_3$ and $A_4$  with the choice $\epsilon_1=\epsilon_2=\epsilon_3=0$, $\epsilon_4=\Delta_4/2$, where $\Delta_4(A_3)/2=0.6572$ and $\Delta_4(A_4)/2=3.9646$. For these matrices, we have that $\sigma_n(A_3)=1.0943$ and  $\sigma_n(A_4)=4.2327$. The bounds obtained applying \eqref{sigma} are $\sigma_n(A_3)>0.3357$ and $\sigma_n(A_4)>0.8680$. 
    \end{example}

  \section{Error bounds for LCP of
Nekrasov matrices}

Given an $n\times n$ real matrix $A$ and $q \in {\bf R}^n$, these problems look for solutions $x^*\in {\bf R}^n$ of
\begin{equation} \label{(1.1)}
Ax+q\ge 0,\quad x\ge 0,\quad x^T(Ax+q)=0. 
\end{equation}
This problem (\ref{(1.1)}) is usually denoted by LCP($A,q$). A real square matrix
is called a {\it $P$-matrix} if all its principal minors are
positive. Let us recall (see \cite{4}) that $A$ is a $P$-matrix if and
only if the LCP($A,q$) (\ref{(1.1)}) has a unique solution $x^*$ for each $q
\in {\bf R}^n$. 

Let $A$ be a real $H$-matrix with all its  diagonal entries
positive. Then $A$ is a $P$-matrix and so  we can apply the third
inequality of Theorem 2.3 of \cite{c-x} and obtain
for any $x\in {\bf R}^n$ the inequality:
$$\|x-x^*\|_\infty \leq {\rm max}_{d\in [0,1]^n} \|(I-D+DA)^{-1}\|_\infty \|r(x)\|_\infty ,$$
where we denote by $I$  the $n \times n$ identity matrix, by $D$ the
diagonal matrix $D={\rm diag}(d_i)_{i=1}^n$ with $0 \leq d_i\leq 1$ for all
$i=1,\ldots,n$,  by $x^*$  the solution of the LCP($A,q$) and by
$r(x):={\rm min}(x,Ax+q)$, where the min operator denotes the
componentwise minimum of two vectors.

By (2.4) of \cite{c-x}, given in Theorem 2.1 of \cite{c-x}, when $A=(a_{ij})_{1\leq i,j\leq n}$ is  a real
$H$-matrix with all its  diagonal entries positive, then we have
 \begin{equation}\label{(2.5)}
 {\rm max}_{d\in [0,1]^n} \|(I-D+DA)^{-1}\|_\infty \leq \| (\mathcal{M}(A))^{-1}{\rm max}(\Lambda,I)\|
_\infty , 
\end{equation}
 where we denote by $\mathcal{M}(A)$  the
comparison matrix of $A$, by $\Lambda$  the diagonal part of $A$
($\Lambda :={\rm diag}(a_{ii})_{i=1}^n$) and by ${\rm max}(\Lambda,I):={\rm
diag} ( {\rm max}\{a_{ii},1\})_{i=1}^n $.

The next theorem, corresponding to Theorem 2.1 of \cite{11}, shows the application of obtaining scaling matrices to transform an $H$-matrix into an SDD matrix in order to derive error bounds for LCP.

\begin{theorem} \label{2.2} Suppose that  $A=(a_{ij})_{1\leq i,j\leq n}$
is an $H$-matrix  with all its  diagonal entries positive. Let
$S= diag( s_i)_{i=1}^n,s_i>0$ for all $i\in N$, be
a diagonal matrix such that $AS$ is SDD. For any $i=1, \ldots,n$, let $\bar
\beta_i:=a_{ii} s_i-\sum_{j\neq i}\mid a_{ij}\mid s_j$.
Then
\begin{equation}\label{(2.6)}
{\rm max}_{d \in [0,1]^n} \|(I-D+DA)^{-1}\|_\infty \leq {\rm max}\left\{{{\rm
max}_i\{ s_i\} \over {\rm min}_i \{  \beta_i \}}, {{\rm
max}_i\{s_i\} \over {\rm min}_i \{s_i \} }\right\}. 
\end{equation}
\end{theorem} 

The following theorem provides an error bound for the
particular LCP associated to a Nekrasov matrix.

\begin{theorem} \label{2.3} Let $A=(a_{ij})_{1 \leq i,j \leq n}$ be a
Nekrasov matrix with all its diagonal entries positive. Let
$S=diag( s_i)_{i=1}^n$ and $\epsilon _i$ ($i\in N$) be the diagonal matrix and  positive real
numbers, respectively, defined in Theorem \ref{corolarioSDD}. Then
\begin{equation}\label{(2.7)}
{\rm max}_{d \in [0,1]^n} \|(I-D+DA)^{-1}\|_\infty \leq {\rm max}\left\{{1 \over {\rm min}_i \{ \epsilon_i-w_i+p_i \}}, {1
\over {\rm min}_i \{s_i \} }\right\}, 
\end{equation}
 where, for each $i\in
N$, $p_i$ and $w_i$ are defined in Theorem \ref{lacota}.
\end{theorem} 

\begin{proof}  Since $A$ is Nekrasov, $s_i<1$  for all $i\in N$ and $A$ is an $H$-matrix. So, we can apply (\ref{(2.6)}) and then  it is sufficient to prove that
$\bar \beta_i=\epsilon_i-w_i+p_i$ for all $i \in N$. For any $i \in N$,
we have
$$\bar\beta_i=a_{ii}{h_i(A)+ \epsilon _i\over a_{ii}}-\sum_{j \in
N\backslash\{i\}}|a_{ij}|{h_j(A)+ \epsilon _j\over
a_{jj}}$$  
and by
(\ref{hs})
 we can write
\begin{align*}
 \bar\beta_i&=\sum_{j =1}^{i-1} |a_{ij}|{h_j(A) \over
a_{jj}}+\sum_{j =i+1}^{n} |a_{ij}| - \sum_{j\in N\backslash
\{i\}}|a_{ij}|{h_j(A)\over a_{jj} }+ \epsilon_i-  \sum_{j\in N\backslash
\{i\}}\epsilon _j{|a_{ij}|\over a_{jj} }=   \\
 &=
\epsilon_i-  \sum_{j =1}^{i-1}{|a_{ij}|\over a_{jj} }+ \sum_{j =i+1}^n|a_{ij}|\left(1-{h_j(A)+\epsilon _j \over
a_{jj}}\right)=\epsilon_i-w_i+p_i. 
\end{align*}
\end{proof}

As a choice of each parameter $\epsilon _i$ ($i\in N$) in Theorem \ref{2.3}, we recommend (as we already did for $\epsilon _n$) to choose the middle point of the interval where it lies (see Theorem \ref{corolarioSDD}). This choice is applied in the following example, where we present a family of matrices for which our bound (\ref{(2.7)}) is a small
constant, in contrast to the bounds of \cite{c-x}, \cite{n1} and \cite{n2}, which can be arbitrarily
large. Observe also that these matrices do not satisfy the necessary hypotheses to apply the bound of \cite{G-E}.
 \begin{example} Let us consider the family of matrices
 \begin{equation*}
 A=\begin{pmatrix}
 K & -K +2 & -1\\ 
 -K & K & 0\\ 
 -K & \frac{-1}{K} & K
 \end{pmatrix},
 \end{equation*}
 where $K>2$. In this case, $h_1(A)=K-1$, $h_2(A)=K-1$ and $h_3(A)=K-1 + \frac{K-1}{K^2}$.
 Then the bound \eqref{(2.5)} (of \cite{c-x}) is equal to $\frac{2K^{3} + 2K}{K^{2} - 1}$ and the bounds of \cite{n1} and \cite{n2} coincide and are equal to $\frac{2K^3+2K^2}{K^2-K+1}$. We can observe that these  bounds are arbitrarily large when $K \to \infty$. However, our new bound remains controlled. In fact, \eqref{(2.7)} with $\epsilon_1=0$, $\epsilon_2=1/2$ and $\epsilon_3=\frac{2K^2 - 2K +3}{4 K^2}$ gives the bound  $\frac{4K^3}{2K^3-2K^2-2K+1} $.
 
  \end{example}

 \section*{Acknowledgements}
 
  This research has been partially supported by  MTM2015-65433-P \linebreak (MINECO/FEDER) Spanish Research Grant and by Gobierno de Arag\'on.

\end{document}